\documentclass[11pt,reqno]{amsart}
\usepackage[T1]{fontenc}
\usepackage[utf8]{inputenc}
\usepackage{lmodern}
\usepackage{amsmath,amssymb,amsthm,mathtools}
\usepackage{graphicx}
\usepackage{float}
\usepackage{tikz}
\usepackage[margin=1in]{geometry}
\usepackage[pagebackref=true,bookmarksopen=true]{hyperref}
\date{\today}

\newtheorem{theorem}{Theorem}
\newtheorem{proposition}[theorem]{Proposition}

\theoremstyle{remark}
\newtheorem{remark}[theorem]{Remark}

\newcommand{\CC}{\mathbb C}
\newcommand{\R}{\mathbb R}
\newcommand{\T}{\mathbb T}

\title{Semigroups with periodic complex-time boundedness and holomorphy regions}
\author{M.A.~Perelmuter}
\address{SCAD Soft Ltd., 3a Osvity Street, Kyiv, 03037, Ukraine}
\email{mikeperelmuter@gmail.com}
\keywords{complex-time semigroups, periodic boundedness regions, holomorphic continuation, Ornstein--Uhlenbeck semigroup, hypoelliptic Laplacian on the circle}
\subjclass[2020]{47D06, 47D03, 35P05, 60J35}

\begin{document}
\maketitle

\begin{abstract}
We study the periodicity of complex-time boundedness and holomorphy domains of operator semigroups. We show that an imaginary-time covariance identity involving an isometry yields periodicity of these domains. The method is applied to the Ornstein--Uhlenbeck operator and the hypoelliptic Laplacian on the circle, providing a spectral explanation of periodicity independent of the explicit geometry of their boundedness regions. We also discuss the distinction between boundedness and contractivity, using the classical heat semigroup as an example.
\end{abstract}

\section{Introduction}

We shall use the terminology \emph{analytic semigroup} and \emph{analyticity} in the standard sectorial sense. A $C_0$-semigroup on a Banach space is called analytic if it admits a holomorphic extension to a sector with vertex at the origin, bounded on every smaller closed subsector. When the complex-time parameter ranges in a more general domain, not necessarily a sector, we shall instead speak of a \emph{holomorphic extension}, a \emph{bounded holomorphic continuation}, and its \emph{bounded-holomorphy region}.

\subsection{The sectorial picture}

There is a classical reason why sectors play a privileged role in the theory of holomorphic semigroups.  A $C_0$-semigroup is initially defined for nonnegative real time, and if one seeks a complex-time extension that itself retains the semigroup structure, it is natural to let the complex parameter range over an additive subsemigroup of $\CC$ containing the positive real half-axis.  Open sectors with vertex at the origin provide the canonical domains of this kind. This point of view is standard and emphasized, for example, in \cite[Section~1.4]{ERS}.

For symmetric diffusion semigroups, Stein's analyticity theorem gives a universal sector of analyticity and contractivity on $L^p$, $1<p<\infty$, with half-angle $\frac{\pi}{2}\left(1-\left|\frac{p-2}{p}\right|\right)$ \cite[Chapter~III, Theorem~1]{Stein}.
Bakry~\cite{Bakry} obtained a larger sector for diffusion semigroups. Related estimates for differential operators were developed by
Fattorini~\cite{FattoriniBook,FattoriniAngle} and Okazawa~\cite{Okazawa}. Liskevich and Perelmuter~\cite{LiskevichPerelmuter} proved that every symmetric submarkovian semigroup admits a contractive holomorphic extension to the strictly larger sector
$$\Sigma_p:=\left\{z\in\CC:\ |\arg z|< \arccos\left|1-\frac{2}{p}\right|\right\}.$$

The sector $\Sigma_p$ is the largest common sector of contractivity for the class of symmetric submarkovian semigroups.  Sharpness is exhibited by examples due to Voigt~\cite{Voigt} and Kunstmann~\cite{Kunstmann2001,Kunstmann2002} (see also the earlier diffusion-semigroup analysis of Bakry~\cite{Bakry}).  Kriegler~\cite{Kriegler} showed that the same optimal sector $\Sigma_p$ remains valid for symmetric contraction semigroups without the positivity assumption (Haase, Kunstmann and Vogt~\cite{HaaseKunstmannVogt} subsequently gave an elementary proof of Kriegler's extension).

The Ornstein--Uhlenbeck semigroup is a sharp example for this sectorial theory.
Chill, Fa\v{s}angov\'a, Metafune and Pallara~\cite{ChillFasangovaMetafunePallara} determined the optimal sector of analyticity for Ornstein--Uhlenbeck semigroups on $L^p$ with respect to the invariant measure, and Maas and van Neerven~\cite{MaasVanNeerven} obtained the corresponding optimal-angle result in infinite dimensions.  In the classical symmetric case this optimal sector is precisely $\Sigma_p$ above.  Thus, if one asks for the largest sector with vertex at the origin on which the semigroup is contractive and holomorphic, the sectorial theory gives a sharp answer.

\subsection{The periodic picture}

The full complex-time domain may nevertheless be larger.  Once a semigroup is known to admit a holomorphic extension to a sector, one may ask for its maximal bounded holomorphic continuation.  There is no a priori reason for this domain to remain a sector.

The complex-time boundedness of the Ornstein--Uhlenbeck semigroup was studied by Weissler~\cite{Weissler}, Epperson~\cite{Epperson}, and Janson~\cite{Janson}. Their results imply that its boundedness region is periodic in the imaginary direction (see also~\cite{vanNeervenPortal}).

For the Ornstein--Uhlenbeck semigroup on $L^p(\gamma_d)$, $1<p<\infty$, where $\gamma_d(x)=(2\pi)^{-d/2}e^{-|x|^2/2}$, for $p\neq2$ let
$$E_p:=\left\{z\in\CC:\ |\sin(\operatorname{Im} z)|\leq\tan\phi_p\,\sinh(\operatorname{Re} z)\right\},\qquad
\phi_p=\arccos|2/p-1|.$$
This is the closed \emph{Epperson region}.  We write
$$E_p^\circ:=\operatorname{int}E_p =\left\{z\in\CC:\ |\sin(\operatorname{Im} z)|<\tan\phi_p\,\sinh(\operatorname{Re} z)\right\}$$
for its interior, which is the bounded-holomorphy domain. For $p=2$, set
$$E_2=\{z\in\CC:\operatorname{Re}z\geq0\},\qquad E_2^\circ=\{z\in\CC:\operatorname{Re}z>0\}. $$
The sector $\Sigma_p$ is the largest sector with vertex at the origin contained in $E_p$. The actual boundedness and contractivity regions extend beyond that sector and are periodic in the imaginary direction.

A different periodic non-sectorial region occurs for the hypoelliptic Laplacian on the circle
$$L_b=\frac12\left(-\frac{\partial^2}{\partial y^2}+y^2-1\right) -by\frac{\partial}{\partial x},\qquad b\neq0,$$
acting on $L^2(\mathbb T\times\mathbb R)$.  Theorem~2.1 of~\cite{MSV} states that $e^{-zL_b}$ is bounded precisely in the $2\pi i$-periodic region
$$\left\{z\in\CC:\ \operatorname{Re}z\geq0,\quad
\cosh(\operatorname{Re} z)-\frac12(\operatorname{Re} z)\,\sinh(\operatorname{Re} z)\leq\cos(\operatorname{Im} z)\right\}.$$
Moreover, at every such $z$ the norm of $e^{-zL_b}$ equals one \cite[(2.6) and Remark~2.2]{MSV}.

\begin{figure}[H]
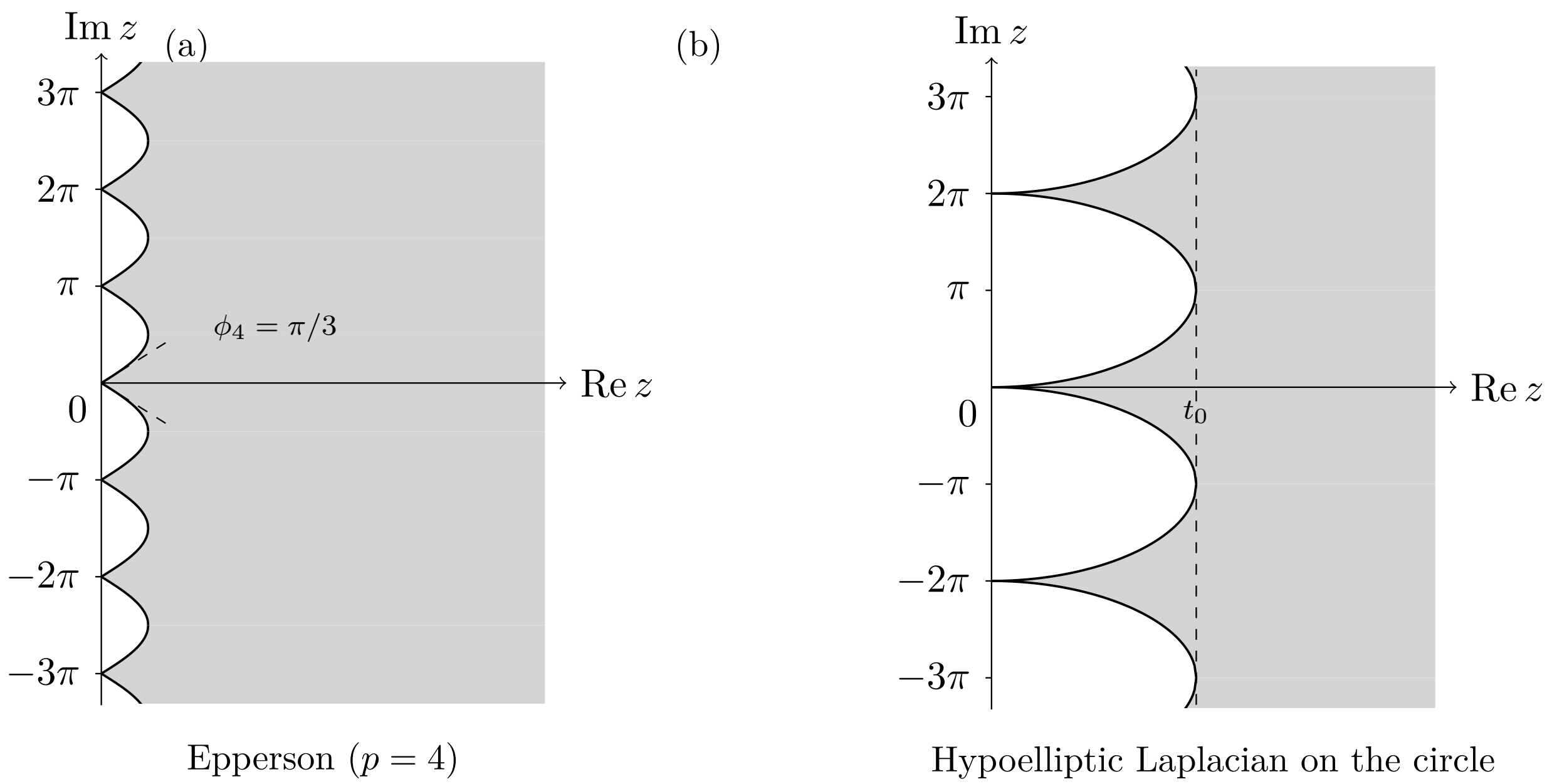

\centering
\resizebox{\textwidth}{!}{%
}
\caption{Periodic complex-time boundedness regions (shaded): (a) the Epperson region for $p=4$; (b) the region for the hypoelliptic Laplacian on the circle. The dashed rays in (a) indicate the optimal sectorial half-angle $\phi_4=\pi/3$. In (b), $t_0\approx2.4$ marks the point beyond which the entire vertical strip to the right is shaded.}
\label{fig:periodic-regions}
\end{figure}

\begin{remark}
The semigroup associated with the hypoelliptic Laplacian on the circle is not analytic in the classical sectorial sense. Indeed, by \cite[Theorem~2.1]{MSV}, for $z=t+is$ with $t>0$ the operator $e^{-zL_b}$ is bounded if and only if
$$\cos s\ge\cosh t-\frac{t}{2}\sinh t.$$
As $t\to0^+$,
$$\cosh t-\frac{t}{2}\sinh t=1-\frac{t^4}{24}+O(t^6),$$
so the boundary near the origin satisfies
$$s=\pm\frac{t^2}{2\sqrt3}+O(t^4).$$
Thus the boundedness region has a quadratic cusp at the origin and contains no sector with vertex at the origin. The semigroup nevertheless has a nontrivial bounded holomorphic continuation on the interior of its boundedness region (see Figure~\ref{fig:periodic-regions}).
\end{remark}

These examples separate two questions that are often conflated in the sectorial theory.  The first is the problem of sectorial analyticity from the positive real axis; the second is the geometry of the
maximal bounded holomorphic continuation.  The explicit results of Epperson and of Mityagin--Siegl--Viola determine the geometry of a fundamental complex-time cell. Our objective is to explain the periodicity of
these boundedness sets independently of their known descriptions, and to distinguish this periodicity from holomorphic continuation. In both examples the repetition is generated by an imaginary-time relation
$$T(z+i\tau)=UT(z),$$
where $U$ is a fixed \emph{isometry}.  Thus the shape of a fundamental cell depends
on the particular operator, whereas its periodic repetition has a simple operator-theoretic explanation.

\section{Periodic boundedness and holomorphic continuation}

$\mathcal L(X)$ denotes the space of bounded linear operators on a complex Banach space $X$.

Our first result isolates a periodicity mechanism --- imaginary-time covariance with an isometry forces periodicity of the boundedness set.

\begin{theorem}
\label{thm:periodic-boundedness}
Let $\mathcal D$ be a dense linear subspace of a complex Banach space $X$, and let $S(z):\mathcal D\to X$, $z\in\CC$, be a family of linear operators.
Suppose that for some $\tau>0$ and a linear isometry $U:X\to X$ one has the \textbf{imaginary-time covariance identity}
$$S(z+i\tau)=US(z),\qquad z\in\CC.$$
Define the boundedness set
$$\mathcal B:=\{z\in\CC:S(z)\text{ has a bounded extension to }X\}.$$
For $z\in\mathcal B$, denote the unique bounded extension by $\widehat S(z)$. Then $\mathcal B+i\tau=\mathcal B$ and
\begin{equation}\label{eq:periodic-boundedness}
\widehat S(z+i\tau)=U\widehat S(z),\qquad
\|\widehat S(z+i\tau)\|=\|\widehat S(z)\|.
\end{equation}
Consequently, for every $M\geq0$, the set $\{z\in\mathcal B:\|\widehat S(z)\|\leq M\}$ is $i\tau$-periodic.
\end{theorem}

\begin{proof}
For $z\in\CC$ and $u\in\mathcal D$, the covariance identity and the isometry of $U$ give
\begin{equation}\label{eq:norm-shift}
\|S(z+i\tau)u\|=\|US(z)u\|=\|S(z)u\|.
\end{equation}
Hence $S(z)$ is bounded on $\mathcal D$ if and only if $S(z+i\tau)$ is, and then
$\sup_{u\in\mathcal D,\,\|u\|\le1}\|S(z+i\tau)u\|=\sup_{u\in\mathcal D,\,\|u\|\le1}\|S(z)u\|$.
In particular, $z\in\mathcal B$ implies $z+i\tau\in\mathcal B$, and $z+i\tau\in\mathcal B$ implies $z\in\mathcal B$. The first implication gives $\mathcal B+i\tau\subset\mathcal B$. Applying the second to $z-i\tau$ gives $\mathcal B\subset\mathcal B+i\tau$. Thus $\mathcal B+i\tau=\mathcal B$.

By the density of $\mathcal D$ in $X$ and the isometry of $U$, both equalities in \eqref{eq:periodic-boundedness} follow.

Finally, by the norm equality \eqref{eq:norm-shift}, $\|\widehat S(z+i\tau)\|\le M$ if and only if $\|\widehat S(z)\|\le M$. Together with $\mathcal B+i\tau=\mathcal B$, this shows that $\{z\in\mathcal B:\|\widehat S(z)\|\le M\}$ is $i\tau$-periodic.
\end{proof}

We next apply this abstract principle to semigroups using a generator-based construction. This construction is used in the examples below.

\begin{proposition}
\label{prop:spectral-realization}
Let $-A$ generate a $C_0$-semigroup $T(t)$ on $X$. Suppose that $X_n\subset \operatorname{Dom}(A)$ are finite-dimensional $A$-invariant subspaces whose span $\mathcal D$ is dense in $X$. Suppose that a
linear isometry $U:X\to X$ satisfies, for some $\tau>0$,
$$U|_{X_n}=e^{-i\tau A_n},\qquad A_n=A|_{X_n}.$$
Then there is a well-defined family $S(z):\mathcal D\to\mathcal D$, $z\in\CC$, such that
$$S(z)|_{X_n}=e^{-zA_n},\qquad S(t)=T(t)|_{\mathcal D}\quad(t\geq0),\qquad S(z+i\tau)=US(z).$$
Moreover, $z\mapsto S(z)u$ is entire for every $u\in\mathcal D$. We also have $[0,\infty)\subset\mathcal B$ and $\widehat S(t)=T(t)$ for $t\ge0$, where $\mathcal B$ and $\widehat S$ are as in Theorem~\ref{thm:periodic-boundedness}.
\end{proposition}

\begin{proof}
Let $Y$ be the sum of finitely many of the spaces $X_n$. Then $Y\subset\operatorname{Dom}(A)$ is finite-dimensional and $A$-invariant. Put $A_Y=A|_Y$. Since each $X_n\subset Y$ is $A_Y$-invariant, $e^{-zA_Y}|_{X_n}=e^{-zA_n}$. If $Y_1$ and $Y_2$ both contain $u$, then $Y_1+Y_2$ is again a finite-dimensional $A$-invariant sum of the spaces $X_n$. Both matrix exponentials agree with $e^{-zA|_{Y_1+Y_2}}$ on their respective subspaces. Hence $S(z)u:=e^{-zA_Y}u$ does not depend on the choice of $Y$ containing $u$. This defines a linear map $S(z):\mathcal D\to\mathcal D$ with $S(z)|_{X_n}=e^{-zA_n}$. For $u\in Y$, the function $z\mapsto e^{-zA_Y}u$ is entire, so $S(\cdot)u$ is entire for every $u\in\mathcal D$.

Let $t\ge0$ and $u\in Y$. The function $v(t)=e^{-tA_Y}u$ is continuously differentiable with values in $Y\subset\operatorname{Dom}(A)$, and satisfies $v'=-Av$, $v(0)=u$. Thus $v$ is a classical solution of the abstract Cauchy problem. Since $T(\cdot)u$ is also a classical solution and such solutions are unique, $T(t)u=v(t)=S(t)u$.

Let $u\in X_n$ and $z\in\CC$. Since $S(z)u\in X_n$,
$$US(z)u=e^{-i\tau A_n}e^{-zA_n}u=e^{-(z+i\tau)A_n}u=S(z+i\tau)u.$$
By linearity, the identity extends to $\mathcal D$. Since $T(t)$ is bounded and agrees with $S(t)$ on the dense set $\mathcal D$, we also have $t\in\mathcal B$ and $\widehat S(t)=T(t)$.
\end{proof}

The preceding results concern boundedness. We now identify the additional condition that yields a periodic holomorphy domain.

\begin{proposition}
\label{prop:local-holomorphy}
In addition to the hypotheses of Theorem~\ref{thm:periodic-boundedness}, assume that $z\mapsto S(z)u$ is entire for every $u\in\mathcal D$. Let $\mathcal H$ be the union of all open subsets
$\Omega\subset\mathcal B$ on which $\widehat S$ is locally bounded in operator norm. Then $\mathcal H$ is open and $i\tau$-periodic, and $z\mapsto\widehat S(z)$ is $\mathcal L(X)$-valued holomorphic on $\mathcal H$. In particular, every open $\Omega\subset\mathcal B$ on which $\|\widehat S(z)\|\leq1$ is contained in $\mathcal H$, so $\widehat S$ is holomorphic and contractive on $\Omega$.
\end{proposition}

\begin{proof}
By definition, $\mathcal H$ is an open subset of $\mathcal B$, and $\widehat S$ is locally bounded on $\mathcal H$. Let $z\in\mathcal H$ and choose an open set $\Omega$ from the defining union such that $z\in\Omega$. By Theorem~\ref{thm:periodic-boundedness}, $\Omega\pm i\tau\subset\mathcal B$ are open and $\|\widehat S(w\pm i\tau)\|=\|\widehat S(w)\|$ for $w\in\Omega$. Hence $\widehat S$ is locally bounded on $\Omega\pm i\tau$, so $\Omega\pm i\tau\subset\mathcal H$. Thus $\mathcal H+i\tau=\mathcal H$. Since $\mathcal H\subset\mathcal B$, identity \eqref{eq:periodic-boundedness} holds on $\mathcal H$.

Let $K\subset\mathcal H$ be compact. Covering $K$ by finitely many neighborhoods on which $\widehat S$ is bounded gives $M_K:=\sup_{z\in K}\|\widehat S(z)\|<\infty$. For $u\in X$ choose $u_j\in\mathcal D$ with $u_j\to u$. Then
$$\sup_{z\in K}\|\widehat S(z)u-S(z)u_j\|\le M_K\|u-u_j\|\to0.$$
So $\widehat S(\cdot)u$ is a locally uniform limit of the entire functions $S(\cdot)u_j$, hence holomorphic on $\mathcal H$ by Weierstrass' theorem on locally uniform limits of holomorphic functions. A locally bounded, strongly holomorphic operator family is holomorphic in operator norm (see, e.g., \cite[Ch.~III, Theorem~3.12]{Kato}).

Finally, if $\Omega\subset\mathcal B$ is open and $\|\widehat S(z)\|\leq1$ on $\Omega$, then $\Omega\subset\mathcal H$, and hence $\widehat S$ is holomorphic and contractive on $\Omega$.
\end{proof}

\section{Examples}

\subsection{The Ornstein--Uhlenbeck operator}

Let
$$\mathrm{d}\gamma_d(x)=(2\pi)^{-d/2}e^{-|x|^2/2}\,\mathrm{d}x$$
and consider on $L^p(\R^d,\gamma_d)$, $1<p<\infty$, the Ornstein--Uhlenbeck operator
$$L=-\Delta+x\cdot\nabla.$$
Let $\mathcal P$ denote the span of the multivariate Hermite polynomials. For each $n\ge0$, let $\mathcal W_n$ denote the $n$-th homogeneous Wiener chaos, which is the span of the multivariate Hermite polynomials of total degree $n$. Then
$$\mathcal P=\bigoplus_{n\geq 0}^{\mathrm{alg}}\mathcal W_n.$$
In $L^2(\gamma_d)$ the spaces $\mathcal W_n$ are mutually orthogonal and the Wiener chaos decomposition (see \cite[Proposition~1.14]{Urbina}) is
$$L^2(\gamma_d)=\bigoplus_{n=0}^{\infty}\mathcal W_n.$$ For $1<p<\infty$, let $-L$ denote the generator of the Ornstein--Uhlenbeck $C_0$-semigroup on $L^p(\gamma_d)$. Each $\mathcal W_n$ lies in $\operatorname{Dom}(L)$, and the operator acts diagonally on these spaces:
$$L|_{\mathcal W_n}=nI,$$
and consequently
$$e^{-zL}|_{\mathcal W_n}=e^{-nz}I.$$

We will obtain the imaginary-time covariance from the reflection. Let
$$Rf(x):=f(-x).$$
Since $\gamma_d$ is invariant under $x\mapsto-x$, $R$ is a surjective isometry of $L^p(\gamma_d)$. The Hermite polynomials in $\mathcal W_n$ have parity $(-1)^n$. Thus $R|_{\mathcal W_n}=(-1)^nI$. Since $L|_{\mathcal W_n}=nI$, we also have
$$e^{-i\pi L}|_{\mathcal W_n}=e^{-i\pi n}I=(-1)^nI=R|_{\mathcal W_n}.$$
Hence $e^{-i\pi L}=R$ on the span $\mathcal P$ of the Hermite polynomials.

\begin{proposition}\label{prop:ou-pi}
Let $1<p<\infty$ and let $S_p(z)$ be the Hermite-polynomial family $S_p(z)|_{\mathcal W_n}=e^{-nz}I$. Its boundedness set
$\mathcal B_p$ and locally bounded holomorphy set $\mathcal H_p$ satisfy
$$\mathcal B_p+i\pi=\mathcal B_p,\qquad \mathcal H_p+i\pi=\mathcal H_p.$$
For every $z\in\mathcal B_p$,
$$\widehat S_p(z+i\pi)=R\widehat S_p(z),\qquad \|\widehat S_p(z+i\pi)\|=\|\widehat S_p(z)\|.$$
Consequently the contractivity set is also $i\pi$-periodic.
\end{proposition}

\begin{proof}
The finite-dimensional Hermite chaoses have dense span in $L^p(\gamma_d)$ (polynomials are dense in $L^p(\gamma_d)$, $1\le p<\infty$) and consist of eigenvectors of the generator with eigenvalue $n$. Reflection $R$ is a surjective isometry and $R|_{\mathcal W_n}=(-1)^nI=e^{-i\pi L|_{\mathcal W_n}}$. Apply Proposition~\ref{prop:spectral-realization}, Theorem~\ref{thm:periodic-boundedness}, and Proposition~\ref{prop:local-holomorphy}.
\end{proof}

By Epperson's characterization of bounded Hermite operators \cite[Theorem~1.1 and Notes~1.2--1.3]{Epperson}, the family $S_p(z)$ with $d=1$ is bounded precisely for $z\in E_p$, and is a contraction there. For $d>1$, $S_p(z)$ is the $d$-fold tensor product of this family. By Fubini's theorem it is a contraction whenever the one-dimensional operator is; conversely, unboundedness in one dimension implies unboundedness in $d$ dimensions by restricting to functions of one variable. Thus $\mathcal B_p=E_p$ for all $d\ge1$. Moreover, $\|\widehat S_p(z)\|=1$ on $E_p$, since $S_p(z)1=1$.

Consequently, by Proposition~\ref{prop:local-holomorphy}, $\widehat S_p$ is holomorphic on $E_p^\circ$ (as an open subset of $\mathcal B_p$ on which $\|\widehat S_p\|\le1$), and $E_p^\circ\subset\mathcal H_p\subset\operatorname{int}\mathcal B_p=E_p^\circ$. Thus $\mathcal H_p=E_p^\circ$.

Thus the periodicity of both the boundedness set and the bounded-holomorphy domain follows directly from the spectral structure of the Ornstein--Uhlenbeck operator, independently of the explicit geometric description of the Epperson region.

\subsection{The hypoelliptic Laplacian on the circle}

Mityagin, Siegl and Viola~\cite{MSV} study the hypoelliptic Laplacian on the circle
$$L_b=\frac12(-\partial_y^2+y^2-1)-by\partial_x,\qquad b\neq0,$$
in the Hilbert space
$$H=L^2(\T_x\times\R_y),\qquad \T=\R/2\pi\mathbb Z.$$
The Fourier decomposition in the circle variable is
$$H=\bigoplus_{n\in\mathbb Z}F_n,\qquad F_n=\{e^{inx}f(y):f\in L^2(\R_y)\},$$
and the spaces $F_n$ are mutually orthogonal and reduce $L_b$ \cite[(1.2), (1.3)]{MSV}.  On $F_n$,
$$L_b|_{F_n}=P_{bn}+\frac12(bn)^2,$$
where $P_a=\frac12\bigl(-\partial_y^2+(y-ia)^2-1\bigr)$ is the shifted harmonic oscillator \cite[(1.4), (1.5)]{MSV}.

For $n\in\mathbb Z$ and $k\in\mathbb N_0$ let $h_k$ be the $k$-th Hermite function \cite[(2.2)]{MSV} and put
$$f_{n,k}(x,y)=e^{inx}h_k(y-ibn)\in F_n.$$
By \cite[(2.1)]{MSV}, these are eigenfunctions of $L_b$:
$$L_b f_{n,k}=\left(k+\frac12b^2n^2\right)f_{n,k}=:\lambda_{n,k}f_{n,k},\qquad n\in\mathbb Z,\ k\in\mathbb N_0.$$
Let $\mathcal D_b=\operatorname{span}\{f_{n,k}:n\in\mathbb Z,\ k\in\mathbb N_0\}$. The functions $f_{n,k}$ are linearly independent: different Fourier modes are orthogonal, and for fixed $n$ the polynomial factors of $h_k(y-ibn)$ have distinct degrees. The eigenfunctions of each $P_a$ form a complete set \cite[Section~1]{MSV}. Together with the orthogonal Fourier decomposition, this shows that $\mathcal D_b$ is dense in $H$. Following \cite[(2.3)]{MSV}, set
$$S_b(z)f_{n,k}=e^{-z\lambda_{n,k}}f_{n,k},\qquad n\in\mathbb Z,\ k\in\mathbb N_0,\ z\in\CC,$$
and extend linearly to $\mathcal D_b$. By \cite[Theorem~2.1]{MSV}, $e^{-zL_b}$ is defined as the graph closure of this family, so $S_b(z)$ has a bounded extension to $H$ if and only if $e^{-zL_b}$ is bounded, and the two extensions coincide.
Define $U_b$ by its action on the orthogonal Fourier decomposition,
$$U_b|_{F_n}=e^{-\pi i b^2n^2}I_{F_n}.$$
Since the multipliers have modulus one, $U_b$ is unitary on $H$.

\begin{proposition}\label{prop:msv-covariance}
Let $S_b$ be the family on $\mathcal D_b$ defined above. Its boundedness set $\mathcal B_b$ and locally bounded holomorphy set $\mathcal H_b$ are $2\pi i$-periodic. Moreover, $\mathcal H_b=\operatorname{int}\mathcal B_b$, and $\widehat S_b$ is holomorphic and contractive on $\mathcal H_b$.
\end{proposition}

\begin{proof}
For every $t>0$,
$$\cosh t-\frac t2\sinh t<1.$$
Hence $t\in\mathcal B_b$ by \cite[Theorem~2.1]{MSV}, and
$\|\widehat S_b(t)\|=1$ by \cite[(2.6)]{MSV}.
The identity $S_b(t+s)=S_b(t)S_b(s)$ on $\mathcal D_b$ therefore extends
to a contraction semigroup $T_b(t)=\widehat S_b(t)$, $T_b(0)=I$.
It is strongly continuous: continuity holds on $\mathcal D_b$ by the
eigenfunction formula and extends to $H$ by density and contractivity.
If $-A$ is its generator, then
$$Af_{n,k}=\lambda_{n,k}f_{n,k}.$$

For $X_{n,k}=\operatorname{span}\{f_{n,k}\}$ and $\tau=2\pi$,
$$U_b|_{X_{n,k}}=e^{-\pi i b^2n^2}I
=e^{-i\tau\lambda_{n,k}}I=e^{-i\tau A|_{X_{n,k}}}.$$
Since $U_b$ is a surjective isometry, Proposition~\ref{prop:spectral-realization}
gives the covariance
$$S_b(z+2\pi i)=U_bS_b(z)\qquad(z\in\CC)$$
on $\mathcal D_b$, as well as the entire dependence of $S_b(z)u$ on $z$
for $u\in\mathcal D_b$.
Theorem~\ref{thm:periodic-boundedness} and
Proposition~\ref{prop:local-holomorphy} now give the periodicity
and holomorphy on $\mathcal H_b$.

Finally, \cite[(2.6) and Remark~2.2]{MSV} imply
$\|\widehat S_b(z)\|=1$ for every $z\in\mathcal B_b$
(the points $z\in2\pi i\mathbb Z$ give unitary operators).
Thus $\widehat S_b$ is locally bounded on
$\operatorname{int}\mathcal B_b$, and
Proposition~\ref{prop:local-holomorphy} yields
$\operatorname{int}\mathcal B_b\subset\mathcal H_b$.
The reverse inclusion follows from the definition of $\mathcal H_b$.
Contractivity on $\mathcal H_b$ follows from the same norm identity.
\end{proof}

The exact boundedness set, including its boundary, is described in~\cite[Theorem~2.1]{MSV}. Proposition~\ref{prop:msv-covariance} establishes periodicity of the boundedness set independently of its geometric description, and shows that its interior is the bounded-holomorphy domain.

\section{Concluding remarks}

Periodicity is not restricted to generators with discrete spectrum.
For $1<p<\infty$, consider
$$X=L^p(\R^d\times(0,\infty),\gamma_d\otimes ds).$$
Let $L=-\Delta_x+x\cdot\nabla_x$ be the Ornstein--Uhlenbeck generator
and let $(Bf)(x,s)=sf(x,s)$. The operator $\widetilde L=L+B$ generates
the semigroup
$$
(e^{-t\widetilde L}f)(x,s)
=e^{-ts}(e^{-tL}f(\cdot,s))(x),\qquad t\ge0.
$$
For complex $z$ with $e^{-zL}$ bounded, the corresponding operator satisfies
$$
\|e^{-z\widetilde L}\|=\sup_{s>0}e^{-s\operatorname{Re}z}\,\|e^{-zL}\|.
$$
Since $E_p\subset\{\operatorname{Re}z\ge0\}$, the boundedness set of
$\widetilde L$ is precisely $E_p$, and its bounded-holomorphy domain
is $E_p^\circ$.
To see periodicity directly, let $\mathcal D$ be the span of functions
$h(x)g(s)$, where $h$ is a Hermite polynomial and
$g\in C_c^\infty(0,\infty)$, and define $S(z)$ on $\mathcal D$ by
$S(z)(h g)=e^{-nz}h\,e^{-zs}g$ when $h\in\mathcal W_n$.
Then
$$
S(z+i\pi)=US(z),\qquad (Uf)(x,s)=e^{-i\pi s}f(-x,s).
$$
The operator $U$ is a surjective isometry of $X$, so
Theorem~\ref{thm:periodic-boundedness} gives periodicity directly.
Moreover, every $\lambda\ge0$ belongs to the approximate point spectrum
of $\widetilde L$, so its spectrum is not discrete.

The exact geometry of a complex-time boundedness region can be investigated using explicit integral kernels, as in the work of Weissler~\cite{Weissler} and Epperson~\cite[Section~2]{Epperson} (Mehler kernel and complex Gaussian operators in dimension one).
Janson~\cite{Janson} gave an alternative proof using stochastic integration and It\^o calculus and established a corresponding criterion for general second-quantization operators.
The method of Mityagin--Siegl--Viola~\cite[Sections~3 and~4]{MSV} suggests a complementary spectral approach, based on the eigenvalues of the generator and the norms of its spectral projections. For a given $t>0$, boundedness of $e^{-(t+is)L_b}$ for all $s\in\R$ is related to absolute convergence of the spectral decomposition of $e^{-tL_b}$ \cite[Theorem~3.6]{MSV}.

This spectral approach concerns the quantitative determination of the shape of the boundedness set, rather than its periodicity. It suggests a general problem for operators with discrete spectrum: to identify spectral or geometric hypotheses under which the distribution of eigenvalues and the growth of the norms of spectral projections determine the precise boundedness region.

Other examples of semigroups with a periodic boundedness set come from quadratic operators on $L^2(\R^n)$. One of them is the operator $Q_\theta=-e^{-i\theta}\partial_x^2+e^{i\theta}x^2$, $|\theta|<\pi/2$, whose boundedness set and norm are computed in~\cite{Viola2016}. (The periodicity of this set is also noted in~\cite[Section~1.3]{AV2018}.) Another class consists of the operators $P_M=\frac12\sum\limits_{j,k=1}^n m_{jk}(-\partial_{x_k}+x_k)(\partial_{x_j}+x_j)$, $M=(m_{jk})$ an $n\times n$ complex matrix, studied in~\cite{AV2015,AV2018}.

A further feature shared by the Ornstein--Uhlenbeck operator and the hypoelliptic Laplacian on the circle is a pointwise dichotomy: at every complex time at which the operator is bounded, it is in fact a contraction. Thus, for these examples, the boundedness and contractivity sets coincide. This property does not hold for general analytic contraction semigroups.
A counterexample is the classical heat semigroup $T(z)=e^{z\Delta}$ on $L^p(\R^d)$, $1<p<\infty$, $p\ne2$. It is a Markov contraction for real $z=t\ge0$ and is bounded for every $\operatorname{Re}z>0$, whereas it is contractive precisely when $|\arg z|\le\phi_p$, where $\phi_p=\arccos|1-2/p|$ (with $z=0$ included). Thus its boundedness region is strictly larger than its contractivity region.

It would be of interest to characterize generators for which this boundedness--contractivity dichotomy holds. To the best of our knowledge, a general structural characterization is not currently available; we regard this as an \emph{open} problem.

\section*{Acknowledgments}
AI tools were used during the preparation of this manuscript to assist with proofreading and checking routine mathematical calculations. The graphics were created using AI tools. The final manuscript was entirely reviewed, revised, and verified by the author, who retains sole responsibility for the final content and mathematical accuracy of this work.


\begin{thebibliography}{99}

\bibitem{AV2015} A.~Aleman, J.~Viola, Singular-value decomposition of solution operators to model evolution equations, \emph{Int. Math. Res. Not. IMRN} (2015), no.~17, 8275--8288.

\bibitem{AV2018} A.~Aleman, J.~Viola, On weak and strong solution operators for evolution equations coming from quadratic operators, \emph{J. Spectr. Theory} \textbf{8} (2018), no.~1, 33--121.

\bibitem{Bakry} D.~Bakry, Sur l'interpolation complexe des semigroupes de diffusion, in \emph{S\'eminaire de Probabilit\'es XXIII}, Lecture Notes in Math., vol.~1372, Springer, Berlin, 1989, 1--20.


\bibitem{ChillFasangovaMetafunePallara} R.~Chill, E.~Fa\v{s}angov\'a, G.~Metafune, D.~Pallara, The sector of analyticity of the Ornstein--Uhlenbeck semigroup on $L^p$ spaces with respect to invariant measure, \emph{J. London Math. Soc.} (2) \textbf{71} (2005), no.~3, 703--722.

\bibitem{ERS} M.~Elin, S.~Reich, D.~Shoikhet, \emph{Numerical Range of Holomorphic Mappings and Applications}, Birkh\"auser, Cham, 2019.

\bibitem{Epperson} J.B.~Epperson, The hypercontractive approach to exactly bounding an operator with complex Gaussian kernel, \emph{J. Funct. Anal.} \textbf{87} (1989), no.~1, 1--30.

\bibitem{FattoriniBook} H.O.~Fattorini, \emph{The Cauchy Problem}, Encyclopedia of Mathematics and its Applications, vol.~18, Addison-Wesley, Reading, MA, 1983; Cambridge University Press, Cambridge, 1984.

\bibitem{FattoriniAngle} H.O.~Fattorini, On the angle of dissipativity of ordinary and partial differential operators, in \emph{Functional Analysis, Holomorphy and Approximation Theory II}, North-Holland Mathematics Studies, vol.~86, North-Holland, Amsterdam, 1984, 85--111.

\bibitem{HaaseKunstmannVogt} M.~Haase, P.C.~Kunstmann, H.~Vogt, On the numerical range of generators of symmetric $L^\infty$-contractive semigroups, \emph{Arch. Math. (Basel)} \textbf{107} (2016), no.~5, 553--559.

%

\bibitem{Janson} S.~Janson, On complex hypercontractivity, \emph{J. Funct. Anal.} \textbf{151} (1997), no.~1, 270--280.

\bibitem{Kato} T.~Kato, \emph{Perturbation Theory for Linear Operators}, 2nd ed., Grundlehren der Mathematischen Wissenschaften, vol.~132, Springer-Verlag, Berlin--New York, 1976.

\bibitem{Kriegler} C.~Kriegler, Analyticity angle for non-commutative diffusion semigroups, \emph{J. London Math. Soc.} (2) \textbf{83} (2011), no.~1, 168--186.

\bibitem{Kunstmann2001} P.C.~Kunstmann, Uniformly elliptic operators with maximal $L^p$-spectrum in planar domains, \emph{Arch. Math. (Basel)} \textbf{76} (2001), no.~5, 377--384.

\bibitem{Kunstmann2002} P.C.~Kunstmann, $L^p$-spectral properties of the Neumann Laplacian on horns, comets and stars, \emph{Math. Z.} \textbf{242} (2002), no.~1, 183--201.

\bibitem{LiskevichPerelmuter} V.A.~Liskevich, M.A.~Perelmuter, Analyticity of submarkovian semigroups, \emph{Proc. Amer. Math. Soc.} \textbf{123} (1995), no.~4, 1097--1104.

\bibitem{MaasVanNeerven} J.~Maas, J.~van Neerven, On analytic Ornstein--Uhlenbeck semigroups in infinite dimensions, \emph{Arch. Math. (Basel)} \textbf{89} (2007), no.~3, 226--236.

\bibitem{MSV} B.~Mityagin, P.~Siegl, J.~Viola, The shifted harmonic oscillator and the hypoelliptic Laplacian on the circle, \emph{Ann. Henri Poincar\'e} \textbf{22} (2021), no.~10, 3313--3355.

\bibitem{Okazawa} N.~Okazawa, Sectorialness of second order elliptic operators in divergence form, \emph{Proc. Amer. Math. Soc.} \textbf{113} (1991), no.~3, 701--706.

\bibitem{Stein} E.M.~Stein, \emph{Topics in Harmonic Analysis Related to the Littlewood--Paley Theory}, Annals of Mathematics Studies, vol.~63, Princeton University Press, Princeton, NJ, 1970.

\bibitem{Urbina} W.~Urbina-Romero, \emph{Gaussian Harmonic Analysis}, Springer Monographs in Mathematics, Springer, Cham, 2019.

\bibitem{vanNeervenPortal} J.~van Neerven, P.~Portal, The Weyl calculus with respect to the Gaussian measure and restricted $L^p$--$L^q$ boundedness of the Ornstein--Uhlenbeck semigroup in complex time, \emph{Bull. Soc. Math. France} \textbf{146} (2018), no.~4, 691--712.

\bibitem{Viola2016} J.~Viola, The norm of the non-self-adjoint harmonic oscillator semigroup, \emph{Integral Equations Operator Theory} \textbf{85} (2016), no.~4, 513--538.

\bibitem{Voigt} J.~Voigt, The sector of holomorphy for symmetric submarkovian semigroups, in \emph{Functional Analysis}, Proceedings of the First International Workshop held at Trier University, Germany, de Gruyter, Berlin, 1996, 449--453.

\bibitem{Weissler} F.B.~Weissler, Two-point inequalities, the Hermite semigroup, and the Gauss--Weierstrass semigroup, \emph{J. Funct. Anal.} \textbf{32} (1979), no.~1, 102--121.

\end{thebibliography}
\end{document}